# Notes on the Polish Algorithm

## Oliver Deiser

**Abstract**

We study, with the help of a computer program, the Polish Algorithm for finite terms satisfying various algebraic laws, e. g., left distributivity a(bc) = (ab)(ac). While the termination of the algorithm for left distributivity remains open in general, we can establish some partial results, which might be useful towards a positive solution. In contrast, we show the divergence of the algorithm for the laws a(bc) = (ab)(cc) and a(bc) = (ab)(a(ac)).

## 1. Introduction

The Polish Algorithm is a canonical way to expand – according to a given law – two terms t and s in order to get expansions t* and s* such that t* is an initial segment of s*, or vice versa. The most interesting law we study here is left distributivity, also known as self-distributivity. The termination of the Polish Algorithm for this law is an important and long-standing open problem.

We use a computer program to get some information on the performance of the Polish Algorithm. In the case of left distributivity, the algorithm terminates unexpectedly fast for most pairs of terms, but there are some exceptional complex pairs.

Our experimental results support the view that the algorithm always terminates for left distributivity and for central duplication, i. e., the law a(bc) = (ab)(bc). On the other hand we can refute the termination of the algorithm for the laws (a(bc)) = (ab)(cc) and a(bc) = (ab)(a(ac)). Moreover, termination for the law a(bc) = (ab)(ca) seems to be unlikely. The negative results show that the termination of the Polish Algorithm is not a mere syntactical feature which holds for a whole family of similar laws.

Concerning left distributivity, we introduce a natural coding of terms, which makes the algorithm quite fast. The effectiveness of this coding is interesting and not well understood.

A superb reference on the subject is Dehornoy's book "Braids and Self-Distributivity" [5]. We give a self-contained definition of the Polish Algorithm below, but the reader will find in this book much more details, further material – including partial convergence results – and an extensive bibliography.

## 2. Terms in right Polish notation

We consider terms in one fixed variable x and one binary multiplication symbol ∗.

**Definition 2.1**  *(terms in right Polish notation)*
   The set **T** of terms (in the variable x) is inductively defined as follows:
   (i) x ∈ **T**,
   (ii) if s, t ∈ **T**, then st∗ ∈ **T**.

Here st∗ is the concatenation of the terms t and s and the symbol ∗. The reason that we use right Polish notation st∗ instead of the more conventional notation (s ∗ t) will become clear when we discuss the Polish Algorithm.
   We also define a multiplication on terms as follows:

s · t = st∗   for s, t ∈ **T**.

   The notions of a *subterm*, *substitution of a term*, … are straightforward. Once again, we ask the reader to consult [5] for details.
   In some of the examples and tables below we will use the symbol 0 to denote multiplication, and 1 to denote the variable, e. g., xx∗x∗ = 11010.

**Definition 2.2**  *(depth of a term)*
   The depth of a term t ∈ **T** is the cardinality of the occurrences of x in t.

Thus depth(x) = 1 and depth(ts∗) = depth(t) + depth(s) for all s, t ∈ T.

**Canonical free structures for some laws on T**

We consider the following laws:

(1)   abc∗∗ = ab∗ac∗∗      left distributivity or ac-law: a(bc) = (ab)(ac)
(2)   abc∗∗ = ab∗bc∗∗      central duplication or bc-law: a(bc) = (ab)(bc)
(3)   abc∗∗ = ab∗ca∗∗      ca-law: a(bc) = (ab)(ca)
(4)   abc∗∗ = ab∗cb∗∗      cb-law: a(bc) = (ab)(cb)
(5)   abc∗∗ = ab∗cc∗∗      right duplication or cc-law: a(bc) = (ab)(cc)
(6)   abc∗∗ = ab∗aac∗∗∗   aac-law: a(bc) = (ab)(a(ac))

Note that that the first five laws are a complete list of all laws of the form a(bc) = (ab)(de), where d, e ∈ { a, b, c }, and d or e is c. The last law is mentioned in [5], and it is included here since it provides an example of divergence. Central duplication is studied extensively in [7].
   Left distributivity implies, e.g., a(b(cd)) = (ab)((ac)(ad)), and obvious generalizations. In right Polish notation, this reads abcd∗∗∗ = ab∗ac∗ad∗∗∗.
   For each law there is a natural notion of distance on **T**:

**Definition 2.3**   *(L-distance $dist_L(u,v)$)*

Let u, v ∈ **T**, n ≥ 0, and let (L) be one of the above laws.
  (i) We say that *v can be reached from v using (L) in n-steps* if there is a sequence $u_0, \ldots, u_n$ in **T** s.t. $u_0 = u$, $u_n = v$ and each $u_{i+1}$ results from $u_i$ by applying the law (L) to a subterm of $u_i$.
  (ii) The *L-distance of u, v, $dist_L(u, v)$*, is the minimal n such that v can be reached from u using (L) in n steps, if such an n exists. Otherwise we let dist(u, v) = ∞.

Applying a law (L) can be expanding or reducing, and the notions of an *(L)-expansion* and *(L)-reduction* are straightforward. For example, abcd∗∗∗e∗ is an 1-step reduction of abc∗bd∗∗∗e∗ for left distributivity, for all a, …, e ∈ T.

We can now define canonical syntactical models for a law (L) on T:

**Definition 2.4**   *(L-equivalence and canonical free monogenic models)*

Let (L) be one of the above laws. We define for u, v ∈ **T**:

u $\sim_L$ v   *iff*   $dist_L(u, v) < \infty$,

$u/\sim_L \cdot v/\sim_L = (u \cdot v)/\sim_L$.

The structure $\langle \mathbf{T}/\sim_L, \cdot \rangle$ is called the *canonical free monogenic model of (L)*.

The structure **W** = $\langle \mathbf{T}/\sim_L, \cdot \rangle$ for L = "left distributivity" is the structure of greatest interest. Free monogenic left distributive structures (i.e., structures isomorphic to **W**) arise naturally in different areas of mathematics: in the theory of large cardinals in set theory, and in the algebraic theory of braids.

Finally, we define a relation on our models:

**Definition 2.5**   *(left factors)*

Let u, v ∈ **T**, and let (L) be one of the above laws. We define:

$u/\sim_L <_L v/\sim_L$   *iff*   $v' = u' u_1 u_2 \ldots u_n$ for some $v' \sim_L v$, $u' \sim_L u$,
           $u_1, \ldots, u_n \in \mathbf{T}$, n ≥ 1,

where, per convention, $u_0 u_1 \ldots u_n$ is $(\ldots(u_0 \cdot u_1) \cdot u_2) \cdot \ldots) \cdot u_n)$.

We then say that *u is a left factor of v modulo the law (L)*.

We suppress the L-indices from now on, whenever L is clear from the context. Note that < is transitive for each law. Moreover, we sometimes also write u < v instead of u/∼ < v/∼.

Consider the two depth 7 terms u = xxx∗x∗xxx∗∗∗∗ and v = xx∗x∗xxxx∗∗∗∗, and L = "left distributivity". Does u < v or v < u or u ∼ v hold? If so, are these possibilities mutually exclusive? The answer is yes. In fact, u ∼ v holds, which is not trivial to see (we used the Polish Algorithm for this).

We now discuss these questions in more detail.

**The word problem**

The *word problem for (L)* can be formulated as follows:

> Is there an algorithm deciding the question:
> Given two terms u and v, is u $\sim_L$ v true?

The word problems for left distributivity and central duplication are solvable. The left factor relation is the key:

**Theorem 2.6**  (*Dehornoy, Laver* [1, 2, 3, 7, 9])
   Let (L) be left distributivity or central duplication. Then $\langle \mathbf{T}/\sim, < \rangle$ is a linear ordering. In particular, the word problem for (L) is solvable.

Comparability and irreflexivity are non-trivial properties of the left factor relation. In case of left distributivity, the first proof of irreflexivity in 1989 used unprovable large cardinal assumptions from set theory. This initiated the study of left distributivity, resulting in a new ordering of braids among other things (see [5] for the interesting history and the role of unprovable set theoretical principles in the first solution of this finite problem).

In order to see how linearity implies the solvability of the word problem, we can use the following algorithm, which appears in [1]:

**Brute force algorithm to solve the word problem for (L)**
   Given u, v ∈ **T** we enumerate all pairs u′, v′ ∈ **T** s.t. u′ $\sim$ u, v′ $\sim$ v. By comparability of <, there will appear a pair u′, v′ s.t.: (1) u′ is a proper initial segment of v′, *or* (2) u′ = v′, *or* (3) v′ is a proper initial segment of u′. In case (2) we stop with result "u and v are equivalent". Otherwise we stop with result "u and v are not equivalent".

We need comparability of $<_L$ (i.e., *at least* one of u < v, u $\sim$ v, v < u is true) to know that the algorithm terminates, and we need irreflexivity of $<_L$ to know that the answer we give is correct.

In [1] irreflexivity was assumed as an algebraic hypothesis. The existence of an irreflexive monogenic left distributive structure was not known until Laver showed that the algebra of elementary embeddings provides an example (or rather *the* example, since any irreflexive monogenic left distributive structure is, in turn, isomorphic to **W**.) Later Dehornoy found a second example given by the operation $s \cdot t = s\ \text{sh}(t)\ \sigma_1\ \text{sh}(s)^{-1}$ on the free braid group generated by $\sigma_1, \sigma_2, \ldots$, where sh is the shift map induced by $\text{sh}(\sigma_n) = \sigma_{n+1}$. In this way irreflexivity of < for left distributivity can be proven by giving the syntactical object **W** a semantical background in different areas of mathematics. One can give a pure syntactical proof of irreflexivity for left distributivity, see [5, VIII], which is important in the study of central duplication.

## 3. The Polish Algorithm

A better algorithm to decide the order-relation < for left distributivity and central duplication is the Polish Algorithm. (Or rather "would be", since termination for this algorithm is open for both laws. For partial convergence results concerning left distributivity see [5; 6].)

The idea of the algorithm is "iterating away the least difference", a method which is applied in set theory with great success. An abstract formulation reads: Given two objects, we search for the first point of disagreement, and then expand the objects in a way that the expansions agree up to and including this point. We can always make progress in this way, but we have to show that we end up with objects comparable w. r. t. to initial segments. In this case the algorithm has reached its goal.

To formulate the algorithm, it is convenient to make some definitions. We view a term $v \in \mathbf{T}$ as a word $v = v(0)\,v(1)\,\ldots\,v(n-1)$, where $n = \text{length}(v)$. For $u, v \in \mathbf{T}$ we set:

$d(u, v)$ = "the least $i < \min(\text{length}(u), \text{length}(v))$ such that $u(i) \neq v(i)$",

if such an i exists, and we let $d(u, v) = \infty$ otherwise. Thus $d = \infty$ iff u is an initial segment of v *or* u = v *or* v is an initial segment of u. Suppose now that $d = d(u, v) < \infty$, and that $v(d) = *$. Then u and v can be uniquely written as follows:

$u = A\,u_0\,u_1\,u_2 * u_3 * \ldots u_n * * B, \quad v = A\,u_0\,u_1 * C,$

where $n \geq 2$, $u_i \in \mathbf{T}$, A, B, C are strings, and the first letter of the subterm $u_2$ is at position d in u.

**Definition 3.1**   *(dynamic range, head and tail at d)*
In this situation, we call $(u_0, \ldots, u_n)$ the *dynamic range of u at position d*, A the *head of u at d*, and B the *tail of u at d*.

We can now define the Polish Algorithm for left distributivity.

**The Polish Algorithm for left distributivity: abc** = ab*ac****

Let $u, v \in \mathbf{T}$, and let $d = d(u,v)$. We stop if $d = \infty$. Thus let $d < \infty$, and assume first that $v(d) = *$. (Thus $u(d) = x$.) Then we set $v' = v$ and

(+)   $u' = A\,u_0\,u_1 * u_0\,u_2 * * u_0\,u_3 * * \ldots u_0\,u_n * * B,$

where $(u_0, \ldots, u_n)$ is the dynamic range of u at d, A is the head of u at d, and B is the tail of u at d. If $u(d) = *$, then we let $u' = u$ and define $v'$ analogous to $u'$ in the first case. A *step in the algorithm* is to go from $(u, v)$ to $(u', v')$. The algorithm, starting on two terms u, v, repeats this steps until, if ever, it reaches terms $u^*, v^*$ such that $d(u^*, v^*) = \infty$. Then $(u^*, v^*)$ is the *result* of the Polish Algorithm, running on $u, v \in \mathbf{T}$.

The Polish Algorithm for the other laws is defined in the same way, where $u'$ in (+) is in each case defined as follows:

**The Polish Algorithm for central duplication: abc∗∗ = ab∗bc∗∗**

(+)  $u' = A\ u_0\ u_1 * u_1\ u_2 * * z_3\ u_3 * * z_4\ u_4 * * \ldots z_n\ u_n * * B$,

where $z_k = u_1\ u_2 * u_3 * u_4 \ldots u_{k-1} *$ for $3 \leq k \leq n$.

**The Polish Algorithm for the law: abc∗∗ = ab∗ca∗∗**

(+)  $u' = A\ u_0\ u_1 * u_2\ u_0 * * u_3\ u_0 * * \ldots u_n\ u_0 * * B$

**The Polish Algorithm for the law: abc∗∗ = ab∗cb∗∗**

(+)  $u' = A\ u_0\ u_1 * u_2\ u_1 * * u_3\ z_3 * * u_4\ z_4 * * \ldots u_n\ z_n * * B$,

where the $z_k$'s are as above.

**The Polish Algorithm for right duplication: abc∗∗ = ab∗cc∗∗**

(+)  $u' = A\ u_0\ u_1 * u_2\ u_2 * * u_3\ u_3 * * \ldots u_n\ u_n * * B$

**The Polish Algorithm for the law: abc∗∗ = ab∗aac∗∗∗**

(+)  $u' = A\ u_0\ u_0\ u_1 * * u_0\ u_0\ u_2 * * * \ldots \ldots u_0\ u_0\ u_n * * * B$

Example: We let $u = 111011000$ and $v = 110110100$. Then the algorithm for left distributivity runs as follows (where "dyn" means "dynamic range"):

Step 0:   111011000  
          110110100           (d = 2, dyn = 1, 1, 1, 110)

Step 1:   1101100111000  
          110110100           (d = 6, dyn = 110, 110, 1)

Step 2:   1101100111000  
          1101100110100       (d = 9, dyn = 1, 1, 1)

Step 3:   110110011011000  
          1101100110100       (d = 11, dyn = 110, 1, 1)

Step 4:   1101100110101101000  
          1101100110100       (d = 12, dyn = 1101100, 11010, 11010)

Step 5:   110110011010011011001101000,  
          1101100110100.

Thus v is a left factor of u, and in particular, u and v are not equivalent w. r. t. to left distributivity.

Clearly, the following holds:

**Lemma 3.2** *(basic properties of the Polish Algorithm)*
Let (u, v), (u′, v′) be a step in the Polish Algorithm for the law (L). Then:
  (i) u′ is an (L)-expansion of u, v′ is an (L)-expansion of v.
  (ii) More exactly, let v′ = v, and let $(u_0, \ldots, u_n)$ be the dynamic range of u at d(u, v), then u′ is an expansion of u by n – 1 applications of (L).
  (iii) d(u′, v′) > d(u, v).

Note that the algorithm is canonical for each law: Let (u, v), (u′, v′) be a step in the algorithm such that u′ is an n-step expansion of u for an n ≥ 1. Then every m-step expansion ū of u which satisfies d(ū, v) > d(u, v) satisfies m ≥ n. Moreover, if m = n, then ū = u.

## 4. Results for left distributivity

The following results where established using a computer program. Our first result supports the hypothesis that the Polish Algorithm terminates for left distributivity:

**Theorem 4.1** *(convergence results for left distributivity)*
The Polish Algorithm for left distributivity terminates for all pairs of depth ≤ 8.

We also proved convergence for all except 60 pairs of depth 9 (there are more than 2000000 pairs of depth ≤ 9).

Table I in the appendix illustrates the performance of the algorithm. For any depth d:
- N(d) is the cardinality of terms of depth ≤ d,
- Σ(d) is the cardinality of pairs of depth ≤ d, where
  "pair" here means a set { s, t } such that s ≠ t. Thus Σ = N(N – 1)/2.
- n(d) is the number of pairs which terminate in exactly n steps.

Example: There are exactly 103 pairs of depth ≤ 8 which terminate in exactly 18 steps. We omit rows where n(d) = 0 for d = 8, 9.

Looking at the statistics we observe that the algorithm terminates very quickly for almost all pairs, but there is more than one local maximum in general. The pattern becomes more striking when we look at the distribution of depth 9 pairs.

Table II lists the most complex pairs for 6 ≤ d ≤ 9. We see an interesting scheme for those "winners" of a class. The table might be helpful in providing a "trivial proof" of termination:

*Possibility of a "trivial proof"*

Find a function f : term-pairs → ordinal-numbers s. t. f(u, v) > f(u′, v′) for each pair of terms u, v, where (u, v), (u′, v′) is a step in the algorithm.

Note that there exists such a function, if the algorithm terminates: Simply let f(u, v) = "number of steps in which the algorithm, running on u, v, terminates".

**Coding**

Table III shows the complexity of some exceptional examples of depth 8 and 9. The depth 9 example ending in terms of length $10^{19}$ cannot be treated without coding. The idea of our coding is: The algorithm is defined as before, but we additionally define a sequence C = ($c_0$, $c_1$, …) of codes during runtime, and we encode the terms according to the codes C at each step.

To give a precise description of the coded version of the algorithm, we need some definitions. First, we replace the multiplication symbol ∗ by 0, and the variable x by 1. Thus our terms are now certain 01-sequences.

**Definition 4.2**  *(codes and encoding)*

A *code* c is a sequence $n_0$ $n_1$ … $n_k$ of natural numbers. If s = $s_0$ is a sequence of natural numbers, and C = ($c_0$, …, $c_k$) is a sequence of codes, then we define by recursion for i ≤ k:

$s_{i+1}$ = the result of the substitution: i for $c_i$ in $s_i$,

i. e., we replace each subsequence of the form $c_i$ in $s_i$ by i (uniformly and simultaneously), where we view i as a sequence of length 1. We now set

ec(s, C) = $s_{k+1}$,

and call ec(s, C) the *encoding of s according to C*.

We will use only code sequences where $c_0$ is 0, $c_1$ is 1, and each $c_k$ is a sequence of numbers < k for each k ≥ 2.

To give an example, we let $s_0$ = 1 1 0 2 0 1 0 1 1 1 0 0 and C = ($c_0$, …, $c_3$) = (0, 1, 1 1 0, 2 2 0). Then:

$s_1$ = $s_2$ = $s_0$,    $s_3$ = 2 2 0 1 0 1 2 0,    $s_4$ = 3 1 0 1 2 0 = ec(s, C).

The coded version of the Polish Algorithm for left distributivity is now defined as follows:

**The coded Polish Algorithm for left distributivity: abc∗∗ = ab∗ac∗∗**

To begin with, set $c_0$ = 0, $c_1$ = 1. Now assume that we have produced terms u, v after n ≥ 0 steps, and that we have produced a code-sequence C = ($c_0$, …, $c_k$) for a k ≥ 1. u and v are now sequences of numbers. We define d = d(u, v) as before as the position of a first difference (or stop the algorithm, if there is none). As before, let u(d), v(d) be the natural numbers of u and v at that position.

*1. case:* $v(d) = 0$. Then we define the dynamic vector $(u_0, \ldots, u_n)$ of u at d as before, where we now treat 0 as the multiplication, and numbers $\geq 1$ as if they were different variables. We then build $u'$ as before. If $u_0$ happens to be a number, we define en(u′, C), en(v, C) as the result of this step. Otherwise we have length($u_0$) > 1. We then introduce a new code $c_{k+1} = u_0$, and define the pair ec(u′, C), ec(v, C) to be the result of the step, where C is now $(c_0, \ldots, c_{k+1})$. Note that $c_{k+1}$ is a sequence of numbers $\leq k$.

*2. case:* $u(d) = 0$. Analogous to case 1.

*3. case:* $0 < u(d) < v(d)$. In this case we replace v(d) at position d in v by $c_{v(d)}$, i.e., we decode the position d in v. Now we repeat the procedure of searching for the first difference, and enter one of the four cases again. (Thus the step is not complete is this case.)

*4. case:* $0 < v(d) < u(d)$. Analogous to case 3.

The rest is as before. Note that the Polish Algorithm, running on u and v, terminates in n steps if and only if the coded version, running on u and v, terminates in n steps. It is for this reason that we do not treat the decoding as a step of the algorithm, but as an operation inside a step.

As in the example given above, let u = 111011000, v = 110110100. The steps for the coded version of the algorithm are:

| | | | |
|---|---|---|---|
| Step 0: | 111011000 | 110110100 | |
| Step 1: | 1101100111000 | 110110100 | |
| Step 2: | 2201200 | 2202100 | with code 2 = 110 |
| Step 3: | 2202200 | 2202100 | |
| Step 4: | 22021021000 | 2202100 | |
| Step 5: | 32100321000 | 32100 | with code 3 = 220. |

As before, the algorithm terminates in 5 steps, and shows that v < u.

We can define a coded version of the Polish Algorithm for the other laws in a similar way, but the details are quite different, since the expansion u′ of u has a different form for the laws (1), (3) and (4). Coding for the law (2) is much the same as for left distributivity, but it turns out to be not as effective as it is for left distributivity.

Our coding is based on the experimental observation that even the most complex examples of depth 8 use just a couple of different $u_0$-terms, where $u_0$ is the first element of the dynamic range. The coding turned out to be quite effective.

We often use the sequence 0, 1, …, 9, a, b, …, z, A, …, Z, … instead of 0, 1, 2, …, 10, 11, … to index codes. This produces a more readable output avoiding spaces or commas to distinguish 1 0 1 1 from 10 11.

The depth 9 example "2:" in table III of the appendix, running for almost 9 million steps, used just 46 different codes: in about 9 million dynamic vectors there are just 45 different $u_0$'s. Looking at the output, one conjectures that the introduction of new codes decreases exponentially.

The coding makes complex computations possible. Moreover, a theoretical understanding of effectiveness of the coding might be the helpful in proving termination.

Table IV in the appendix shows a typical result of a more complex coded computation.

## 5.  Results for the other laws

As expected, central duplication turns out to be considerably simpler than left distributivity:

**Theorem 5.1**  *(convergence result for central duplication)*
The Polish Algorithm for central duplication abc∗∗ = ab∗bc∗∗ terminates for all pairs of depth ≤ 10.

Table V illustrated the performance of the algorithm. It is interesting to compare it with left distributivity.

**Theorem 5.2**  *(results for the ca-law; exceptional depth 7 pair)*
The Polish Algorithm for the law abc∗∗ = ab∗ca∗ terminates for all pairs of depth ≤ 6. It terminates for all but one pair of depth ≤ 7 in at most 19 steps. The exception is the pair 1111101010000, 1101011101000, which does not terminate in 6666 steps.

After 6666 steps, the algorithm has produced terms of length 5820919 and 20309879. The agreement between these terms is only 63821 or about 1%. We conjecture that the algorithm diverges for this pair.

**Theorem 5.3**  *(convergence result for the cb-law)*
The Polish Algorithm for the law abc∗∗ = ab∗cb∗ terminates for all pairs of depth ≤ 8.

The performance of the algorithm for this law is very similar to that for left distributivity.

Concerning right duplication we could establish divergence:

**Theorem 5.4a** *(divergence for right duplication)*
The Polish Algorithm for the law abc∗∗ = ab∗cc∗∗ terminates for all pairs of depth ≤ 6. It diverges for some pairs of depth ≤ 7.

To illustrate this result, we prove divergence for the simplest pair of terms:

**Theorem 5.4b** *(first pair of divergent terms for right duplication)*
The Polish Algorithm for the law abc∗∗ = ab∗cc∗∗ diverges for the pair
$u_0$ = 11111001000, $v_0$ = 1101011110000.

**Proof**
We set:

a  =  11010110,
b  =  110010111001000,
c  =  1100011011001101110000,
d  =  11000110011100010111001000,
z  =  110001100110.

Let $u_n$, $v_n$ be the terms produced by the algorithm after step n. Then an induction shows that:

$u_{3m+1}$ = a $z^{m-1}$ b,
$v_{3m+1}$ = a $z^{m-1}$ c,
$u_{3m+2}$ = a $z^{m-1}$ d,
$v_{3m+2}$ = a $z^{m-1}$ c,
$u_{3m+3}$ = a $z^{m-1}$ d,
$v_{3m+3}$ = a $z^m$ c,

for all m ≥ 2, where $z^n$ is the concatenation of n copies of z. Thus the algorithm diverges for this pair.

In fact, the algorithm terminates for all but nine pairs of depth ≤ 7 in at most 20 steps. Table VI lists those pairs. We can prove divergence for these and other, more complicated pairs in a similar fashion, introducing appropriate codes like a, b, c, … in the above proof.

Finally, we have the following divergence result for law (6):

**Theorem 5.5** *(divergence result for the aac-law)*
The Polish Algorithm for the law abc∗∗ = ab∗aac∗∗∗ converges for all pairs of depth ≤ 5. It diverges for some pairs of depth ≤ 6, e.g., for the pair 1110100, 11010111000.

The proof of divergence is similar to the proof for right duplication, but a more complicated coding sequence is used.

## Table I    Performance of the Polish Algorithm for left distributivity

| d | 3 | 4 | 5 | 6 | 7 | 8 |
|---|---|---|---|---|---|---|
| N | 3 | 8 | 22 | 64 | 197 | 625 |
| Σ | 3 | 28 | 231 | 2016 | 19110 | 195000 |
| 0 | 1 | 5 | 19 | 67 | 232 | 806 |
| 1 | 1 | 9 | 45 | 205 | 845 | 3399 |
| 2 | 1 | 7 | 46 | 271 | 1335 | 6430 |
| 3 |   | 7 | 62 | 388 | 2295 | 12387 |
| 4 |   |   | 55 | 533 | 3697 | 22771 |
| 5 |   |   | 4 | 454 | 4678 | 35279 |
| 6 |   |   |   | 63 | 4101 | 43429 |
| 7 |   |   |   | 21 | 939 | 38711 |
| 8 |   |   |   | 11 | 373 | 11668 |
| 9 |   |   |   | 2 | 318 | 6501 |
| 10 |   |   |   | 1 | 122 | 5071 |
| 11 |   |   |   |   | 65 | 2634 |
| 12 |   |   |   |   | 73 | 2003 |
| 13 |   |   |   |   | 18 | 1303 |
| 14 |   |   |   |   | 6 | 986 |
| 15 |   |   |   |   | 10 | 441 |
| 16 |   |   |   |   |   | 229 |
| 17 |   |   |   |   |   | 248 |
| 18 |   |   |   |   | 1 | 103 |
| 19 |   |   |   |   | 1 | 141 |
| 20 |   |   |   |   |   | 76 |
| 21 |   |   |   |   | 1 | 154 |
| 22 |   |   |   |   |   | 46 |
| 23 |   |   |   |   |   | 31 |
| 24 |   |   |   |   |   | 24 |
| 25 |   |   |   |   |   | 10 |
| 26 |   |   |   |   |   | 29 |
| 27 |   |   |   |   |   | 11 |
| 28 |   |   |   |   |   | 4 |
| 29 |   |   |   |   |   | 33 |
| 30 |   |   |   |   |   | 8 |

|  | 8 |
|---|---|
| 34 | 1 |
| 36 | 1 |
| 38 | 1 |
| 39 | 1 |
| 41 | 1 |
| 42 | 1 |
| 43 | 2 |
| 44 | 2 |
| 46 | 1 |
| 48 | 2 |
| 50 | 7 |
| 51 | 1 |
| 52 | 6 |
| 56 | 1 |
| 58 | 1 |
| 59 | 1 |
| 60 | 1 |
| 61 | 2 |
| 78 | 1 |
| 473 | 1 |
| 1831 | 1 |

|   | 9 |
|---|---|
| **N** | 2055 |
| Σ | 2110485 |
| 0 | 2877 |
| 1 | 13380 |
| 2 | 28959 |
| 3 | 62938 |
| 4 | 126490 |
| 5 | 234457 |
| 6 | 350388 |
| 7 | 421959 |
| 8 | 379079 |
| 9 | 141392 |
| 10 | 93059 |
| 11 | 73342 |
| 12 | 43394 |
| 13 | 35538 |
| 14 | 26656 |
| 15 | 18208 |
| 16 | 13785 |
| 17 | 8960 |
| 18 | 6089 |
| 19 | 5787 |
| 20 | 4049 |
| 21 | 3915 |
| 22 | 2695 |
| 23 | 2346 |
| 24 | 1211 |
| 25 | 1018 |
| 26 | 693 |
| 27 | 733 |
| 28 | 1018 |
| 29 | 589 |
| 30 | 768 |

| 31 | 488 |
|---|---|
| 32 | 166 |
| 33 | 120 |
| 34 | 103 |
| 35 | 118 |
| 36 | 69 |
| 37 | 39 |
| 38 | 63 |
| 39 | 44 |
| 40 | 76 |
| 41 | 27 |
| 42 | 52 |
| 43 | 56 |
| 44 | 94 |
| 45 | 53 |
| 46 | 38 |
| 47 | 42 |
| 48 | 51 |
| 49 | 74 |
| 50 | 75 |
| 51 | 112 |
| 52 | 148 |
| 53 | 158 |
| 54 | 85 |
| 55 | 37 |
| 56 | 159 |
| 57 | 56 |
| 58 | 43 |
| 59 | 31 |
| 60 | 45 |
| 61 | 52 |
| 62 | 105 |
| 63 | 274 |
| 64 | 20 |
| 65 | 113 |

| 66 | 406 |
|---|---|
| 67 | 15 |
| 68 | 97 |
| 69 | 9 |
| 70 | 3 |
| 71 | 3 |
| 72 | 1 |
| 73 | 1 |
| 76 | 4 |
| 77 | 15 |
| 78 | 12 |
| 79 | 9 |
| 80 | 9 |
| 81 | 5 |
| 82 | 3 |
| 83 | 6 |
| 84 | 7 |
| 85 | 1 |
| 86 | 2 |
| 87 | 26 |
| 88 | 1 |
| 89 | 1 |
| 90 | 2 |
| 91 | 8 |
| 92 | 1 |
| 96 | 1 |
| 97 | 2 |
| 98 | 2 |
| 99 | 2 |
| 100 | 1 |
| 103 | 1 |
| 104 | 1 |
| 105 | 1 |
| 111 | 2 |
| 113 | 2 |

| 114 | 6 |
|---|---|
| 118 | 1 |
| 124 | 1 |
| 125 | 7 |
| 141 | 1 |
| 143 | 1 |
| 146 | 1 |
| 149 | 1 |
| 155 | 1 |
| 160 | 1 |
| 163 | 1 |
| 165 | 6 |
| 166 | 2 |
| 175 | 1 |
| 182 | 1 |
| 183 | 7 |
| 184 | 1 |
| 185 | 33 |
| 186 | 7 |
| 187 | 1 |
| 196 | 1 |
| 197 | 2 |
| 198 | 2 |
| 199 | 6 |
| 200 | 1 |
| 202 | 1 |
| 203 | 1 |
| 204 | 1 |
| 205 | 2 |
| 207 | 1 |
| 208 | 1 |
| 209 | 1 |
| 210 | 2 |
| 217 | 22 |
| 218 | 19 |

| | | | | | |
|---|---|---|---|---|---|
| 219 | 100 | 495 | 4 | 1043 | 3 |
| 220 | 27 | 496 | 4 | 1060 | 1 |
| 221 | 18 | 497 | 3 | 1072 | 1 |
| 222 | 128 | 519 | 2 | 1074 | 1 |
| 223 | 3 | 520 | 2 | 1093 | 1 |
| 224 | 27 | 521 | 1 | 1094 | 7 |
| 225 | 5 | 590 | 1 | 1370 | 1 |
| 226 | 6 | 591 | 6 | 1797 | 1 |
| 227 | 6 | 594 | 1 | 1799 | 6 |
| 228 | 7 | 595 | 7 | 1800 | 2 |
| 229 | 2 | 597 | 1 | 1828 | 1 |
| 230 | 3 | 602 | 1 | 1829 | 8 |
| 231 | 1 | 605 | 1 | 1830 | 9 |
| 235 | 1 | 658 | 1 | 1831 | 41 |
| 236 | 1 | 659 | 7 | 1832 | 17 |
| 237 | 1 | 661 | 1 | 1833 | 14 |
| 238 | 1 | 662 | 7 | 1834 | 1 |
| 239 | 2 | 707 | 1 | 1853 | 1 |
| 301 | 1 | 804 | 1 | 1854 | 2 |
| 302 | 6 | 806 | 1 | 1855 | 2 |
| 305 | 1 | 861 | 1 | 1856 | 6 |
| 306 | 7 | 863 | 1 | 1857 | 1 |
| 312 | 1 | 1008 | 1 | 2738 | 1 |
| 316 | 1 | 1009 | 1 | 2739 | 1 |
| 320 | 1 | 1010 | 3 | 2740 | 1 |
| 463 | 1 | 1011 | 3 | 2741 | 2 |
| 464 | 2 | 1012 | 1 | 4159 | 1 |
| 471 | 1 | 1013 | 2 | 4160 | 1 |
| 472 | 8 | 1014 | 1 | 5070 | 1 |
| 473 | 8 | 1036 | 1 | 7588 | 1 |
| 474 | 10 | 1038 | 1 | ≥ 10000 | 71 |
| 475 | 14 | 1039 | 2 | | |
| 477 | 1 | 1040 | 4 | | |
| 481 | 1 | 1041 | 4 | | |
| 494 | 1 | 1042 | 2 | | |

## Table II   Complex terms for left distributivity

**Depth 6**
1:   11101110000   11010111000   10 steps

**Depth 7**
1:   1110101110000   1101010111000   21 steps

**Depth 8**
1:   111010101110000   110101011110000    473 steps
2:   111010101110000   110101010111000    1831 steps

**Depth 9**
The 71 depth ≤ 9 pairs which do not terminate in 10000 steps are:

```
 1:  111010101110000    11010101111100000     8977140 steps
 2:  111010101110000    11010101111010000     8977136 steps
 3:  111010101110000    11010101011110000     1615928 steps
 4:  111010101110000    11010101010111000
 5:  110101011110000    11101010101110000
 6:  110101010111000    11101010101110000
 7:  110101010111000    11100101011110000     17129 steps
 8:  110101010111000    11100101011101000     17127 steps
 9:  1110101011110000   11010101111100000
10:  1110101011110000   11010101111010000
11:  1110101011110000   11010101011110000
12:  1110101011110000   11010101010111000
13:  1110101011010000   11010101111100000
14:  1110101011010000   11010101111010000
15:  1110101011010000   11010101011110000
16:  1110101011010000   11010101010111000
17:  1110101011010000   11010101011110000
18:  1110101011010000   11010101010111000
19:  1110101011001000   11010101111100000
20:  1110101011001000   11010101111010000
21:  1110101011001000   11010101011110000
22:  1110101011001000   11010101010111000
23:  1110101100111000   11010101111100000     11216 steps
24:  1110101100111000   11010101111010000     11212 steps
25:  1110101010111000   11010101111100000
26:  1110101010111000   11010101111010000
27:  1110101010111000   11010101110110000
28:  1110101010111000   11010101111001000
29:  1110101010111000   11010101101110000
30:  1110101010111000   11010101111000100
```

| | | |
|---|---|---|
| 31: | 1110101010111000  11010101011110000 | |
| 32: | 11101010101110000  11010101011101000 | |
| 33: | 11101010101110000  11010101011011000 | |
| 34: | 11101010101110000  11010101011100100 | |
| 35: | 11101010101110000  11010101100111000 | |
| 36: | 11101010101110000  11010110010111000 | |
| 37: | 11101010101110000  11010101010111000 | |
| 38: | 11101010101110000  11010101111000010 | |
| 39: | 11101010101110000  11010101011100010 | |
| 40: | 11101010111000100  11010101111100000 | |
| 41: | 11101010111000100  11010101111010000 | |
| 42: | 11101010111000100  11010101011110000 | |
| 43: | 11101010111000100  11010101010111000 | |
| 44: | 11011010101110000  11010101011110000 | 57414 steps |
| 45: | 11011010101110000  11010101010111000 | |
| 46: | 11010101111100000  11101010011110000 | |
| 47: | 11010101111100000  11101001011110000 | |
| 48: | 11010101111100000  11100101011110000 | |
| 49: | 11010101111100000  11101010111000010 | |
| 50: | 11010101111010000  11101010011110000 | |
| 51: | 11010101111010000  11101001011110000 | |
| 52: | 11010101111010000  11100101011110000 | |
| 53: | 11010101111010000  11101010111000010 | |
| 54: | 11101010011110000  11010101011110000 | |
| 55: | 11101010011110000  11010101010111000 | |
| 56: | 11101001011110000  11010101011110000 | |
| 57: | 11101001011110000  11010101010111000 | |
| 58: | 11100101011110000  11010101011110000 | |
| 59: | 11100101011110000  11010101011101000 | 17131 steps |
| 60: | 11100101011110000  11010101011011000 | 17128 steps |
| 61: | 11100101011110000  11010101011100100 | 17131 steps |
| 62: | 11100101011110000  11010101010111000 | |
| 63: | 11100101011110000  11010101011100010 | 17130 steps |
| 64: | 11100101011101000  11010101011110000 | |
| 65: | 11100101011101000  11010101011101000 | 17129 steps |
| 66: | 11100101011101000  11010101011011000 | 17126 steps |
| 67: | 11100101011101000  11010101011100100 | 17129 steps |
| 68: | 11100101011101000  11010101010111000 | |
| 69: | 11100101011101000  11010101011100010 | 17128 steps |
| 70: | 11010101011110000  11101010111000010 | |
| 71: | 11010101010111000  11101010111000010 | |

The pairs without an entry in the third column do not terminate in 100000 steps.

## Table III   Some complex examples

**Depth 8**
1:  L: 111010101110000   R: 110101011110000,   terminates in 473 steps
Length of result:
   Left side:    72823933
   Right side:   72823685
Thus R is a left factor of L.

2:  L:  111010101110000   R: 110101010111000,   terminates in 1831 steps
Length of result:
   Left side:    13728381775
   Right side:        87441947
Thus R is a left factor of L.

**Depth 9**
2:  L: 111010101110000   R: 11010101111010000,
                                                   terminates in 8977136 steps

Length of result:
Both sides have length $\sim 1.023 \cdot 10^{19}$.
In fact, R is a left factor of L.

## Table IV   The result of a coded computation for left distributivity

This is a typical coded computation for left distributivity.

We compare the terms 111010101110000 and 110101011110000,
using the law abc** = ab*ac**. The Polish Algorithm terminates after 473 steps.

Codes used:

| | | | | |
|---|---|---|---|---|
| 2 | := | 110 | (code number 2) | Real length: 3 |
| 3 | := | 210 | (code number 3) | Real length: 5 |
| 4 | := | 310 | (code number 4) | Real length: 7 |
| 5 | := | 440 | (code number 5) | Real length: 15 |
| 6 | := | 540 | (code number 6) | Real length: 23 |
| 7 | := | 620 | (code number 7) | Real length: 27 |
| 8 | := | 720 | (code number 8) | Real length: 31 |
| 9 | := | 87100 | (code number 9) | Real length: 61 |
| a | := | 97100971000971000 | (code number 10) | Real length: 275 |
| b | := | a7100 | (code number 11) | Real length: 305 |
| c | := | b6100b61000 | (code number 12) | Real length: 663 |
| d | := | cb200cb1000cb1000cb200cb1000cb10000cb200cb1000cb10000cb1000 | | |
| | | | (code number 13) | Real length: 9725 |
| e | := | db0 | (code number 14) | Real length: 10031 |
| f | := | e90 | (code number 15) | Real length: 10093 |
| g | := | fe71000fe710000fe710000 | (code number 16) | Real length: 60467 |
| h | := | ge71000 | (code number 17) | Real length: 70529 |
| i | := | he61000he610000he2000he61000he610000he10000he61000he610000~ he10000he61000he610000he2000he61000he610000 he10000he61000~ he610000he100000he61000he610000he2000he61000he610000he10000~ he61000he610000he100000he61000he610000he10000 | | |
| | | | (code number 18) | Real length: 2417405 |
| j | := | ie61000ie610000 | (code number 19) | Real length: 4854927 |

Resulting terms:
Left:
   jie2000jie10000jie10000jie2000jie10000jie100000jie2000jie10000jie100000~
   jie100009086100086400086400086400864000
Right:
   jie2000jie10000jie10000jie2000jie10000jie100000jie2000jie10000jie100000jie10000

Real length of resulting terms:
   Left:   72823933
   Right: 72823685

## Table V
## Performance of the Polish Algorithm for central duplication

| d  | 3 | 4  | 5   | 6    | 7     | 8      | 9       | 10       |
|----|---|----|-----|------|-------|--------|---------|----------|
| N  | 3 | 8  | 22  | 64   | 197   | 625    | 2055    | 6917     |
| Σ  | 3 | 28 | 231 | 2016 | 19110 | 195000 | 2110485 | 23918986 |
| 0  | 1 | 5  | 19  | 67   | 232   | 806    | 2806    | 9878     |
| 1  | 1 | 9  | 45  | 197  | 812   | 3228   | 12554   | 48692    |
| 2  | 1 | 7  | 57  | 335  | 1700  | 7882   | 35522   | 154058   |
| 3  |   | 7  | 54  | 450  | 2887  | 16503  | 85086   | 416646   |
| 4  |   |    | 52  | 459  | 4012  | 27822  | 172595  | 975204   |
| 5  |   |    | 2   | 404  | 4077  | 37765  | 279927  | 1866877  |
| 6  |   |    | 1   | 54   | 3493  | 37709  | 367035  | 2894635  |
| 7  |   |    | 1   | 28   | 845   | 32400  | 368219  | 3682822  |
| 8  |   |    |     | 14   | 549   | 11864  | 315910  | 3715190  |
| 9  |   |    |     | 7    | 262   | 8052   | 151903  | 3239598  |
| 10 |   |    |     | 1    | 143   | 4767   | 111227  | 1864254  |
| 11 |   |    |     |      | 61    | 2899   | 73479   | 1442422  |
| 12 |   |    |     |      | 22    | 1495   | 49157   | 1044115  |
| 13 |   |    |     |      | 9     | 846    | 29470   | 750555   |
| 14 |   |    |     |      | 3     | 442    | 20535   | 513043   |
| 15 |   |    |     |      | 2     | 226    | 12498   | 372568   |
| 16 |   |    |     |      | 1     | 129    | 7813    | 259893   |
| 17 |   |    |     |      |       | 80     | 5337    | 187824   |
| 18 |   |    |     |      |       | 36     | 3379    | 134634   |
| 19 |   |    |     |      |       | 25     | 2059    | 96629    |
| 20 |   |    |     |      |       | 15     | 1400    | 67275    |
| 21 |   |    |     |      |       | 5      | 951     | 50860    |
| 22 |   |    |     |      |       | 3      | 550     | 35341    |
| 23 |   |    |     |      |       | 1      | 327     | 25201    |
| 24 |   |    |     |      |       | 1      | 219     | 17780    |
| 25 |   |    |     |      |       | 1      | 169     | 13647    |
| 26 |   |    |     |      |       |        | 124     | 9727     |
| 27 |   |    |     |      |       |        | 55      | 7030     |
| 28 |   |    |     |      |       |        | 61      | 4952     |
| 29 |   |    |     |      |       |        | 40      | 4087     |
| 30 |   |    |     |      |       |        | 25      | 3063     |

| | 9 | 10 |
|---|---|---|
| 31 | 13 | 2256 |
| 32 | 7 | 1623 |
| 33 | 8 | 1185 |
| 34 | 7 | 1066 |
| 35 | 5 | 830 |
| 36 | 5 | 680 |
| 37 | 2 | 611 |
| 38 | 4 | 406 |
| 39 | 2 | 383 |
| 40 | | 320 |
| 41 | | 197 |
| 42 | | 152 |
| 43 | | 112 |
| 44 | | 116 |
| 45 | | 93 |

| | 10 |
|---|---|
| 46 | 71 |
| 47 | 49 |
| 48 | 62 |
| 49 | 54 |
| 50 | 30 |
| 51 | 34 |
| 52 | 34 |
| 53 | 14 |
| 54 | 16 |
| 55 | 24 |
| 56 | 13 |
| 57 | 4 |
| 58 | 5 |
| 59 | 8 |
| 60 | 5 |

| | 10 |
|---|---|
| 61 | 2 |
| 62 | 3 |
| 63 | 4 |
| 64 | 5 |
| 65 | 3 |
| 66 | 2 |
| 67 | 5 |
| 68 | 3 |
| 69 | 1 |
| 70 | 1 |
| 72 | 1 |
| 73 | 1 |
| 75 | 1 |
| 76 | 1 |

**Table VI   Divergence for right duplication**

The Polish Algorithm for right duplication diverges for the following nine pairs of depth ≤ 7:

1:  11111001000   1101011110000  
2:  1111101001000   1101011110000  
3:  1110011001000   1101011110000  
4:  1111100101000   1101011110000  
5:  1110111001000   1101110110000  
6:  1110101101000   1101111010000  
7:  1110101101000   1101101101000  
8:  1111100100100   1101011110000  
9:  1101011110000   1111100100010